\documentclass{gtart}     
\usepackage{amssymb,amsmath,amsthm,psfig,epsfig,graphicx,color}

\newtheorem{thm}{Theorem}
\newtheorem{prop}[thm]{Proposition}
\newtheorem{lemma}[thm]{Lemma}
\newtheorem{cor}[thm]{Corollary}

\theoremstyle{definition}
\newtheorem{defn}[thm]{Definition}

\newtheorem{question}[thm]{Question}

\def\Cal#1{{\cal#1}}
\def\<{\langle}\def\>{\rangle}
\def\what{\widehat}
\def\Z{{\mathbb Z}}
\def\N{{\mathbb N}} 
\def\R{{\mathbb R}}

\def\H{{\mathbb H}}

\def\ga{\gamma}			\def\Ga{\Gamma}

\def\ep{\epsilon}               
\def\Ups{\Upsilon}               
\def \co{\colon\thinspace}
\def \P{\mathcal{P}}
\def \DiffD{\textsl{Diff}(D^2,\partial D^2,\textsl{vol})}
\def \Diff{\textsl{Diff}}
\def \Homeo{\textsl{Homeo}(D^2,\partial D^2,\textsl{vol})}
\def \Mod{\textsl{Mod}}\def \PMod{\textsl{PMod}}
\def\comp{\textsl{complexity}}

\begin{document}

\title{Quasi-isometrically embedded subgroups of braid and diffeomorphism 
groups}
\shorttitle{Quasi-isometrically embedded subgroups of braid groups}

\authors{John Crisp\\Bert Wiest}
\address{Institut de Math\'emathiques de Bourgogne (IMB), 
UMR 5584 du CNRS,\\ 
Universit\'e de Bourgogne, 9 avenue Alain Savary, B.P. 47870,\\ 
21078 Dijon cedex, France\\
\smallskip
IRMAR, UMR 6625 du CNRS,\\ Campus de Beaulieu, Universit\'e de Rennes 1,\\
35042 Rennes, France}
\email{jcrisp@u-bourgogne.fr\\ bertold.wiest@math.univ-rennes1.fr}

\begin{abstract}
We show that  a large class of right-angled Artin groups
(in particular, those with planar complementary defining graph) can be 
embedded quasi-isometrically in pure braid groups and in the group 
$\DiffD$ of area preserving diffeomorphisms of the disk fixing the 
boundary (with respect to the $L^2$-norm metric); this extends results of 
Benaim and Gambaudo who gave quasi-isometric 
embeddings of $F_n$ and $\Z^n$ for all $n>0$. As a consequence we 
are also able to embed a variety of Gromov hyperbolic groups 
quasi-isometrically in pure braid groups and in the 
group $\DiffD$. Examples include hyperbolic surface groups, some 
HNN-extensions of these along cyclic subgroups and the fundamental 
group of a certain closed hyperbolic 3-manifold.
\end{abstract}


\makeshorttitle    


\section{Introduction}\label{S:intro}


We recall that a right-angled Artin group is a group which can be
described by a presentation with a finite number of generators, and
a finite list of relations, each of which states that some pair of
generators commutes. Thus free groups and free abelian groups are
examples of right-angled Artin groups. To any simplicial graph $\Delta$
we can associate a right-angled Artin group
$G(\Delta)$ by having generators corresponding to the vertices of
$\Delta$, and a commutation relation between two generators if and
only if the corresponding vertices of $\Gamma$ are connected by an edge.
Thus, if $\Delta$ has vertex set $\{ 1,2,..,n\}$ then
\[
G(\Delta) =\< a_1,a_2,..,a_n\mid a_ia_j=a_ja_i \text{ for each edge } \{ i,j\} 
\text{ of } \Delta \>\,.
\]

Let $\Ga,\Ga'$ denote groups which are equipped with left-invariant 
metrics $d,d'$ respectively. 
By a {\it quasi-isometric embedding} of $(\Ga,d)$ into $(\Ga',d')$ 
we shall mean a group homomorphism $\phi\co\Ga\to\Ga'$ which is injective 
and which, for some uniform constants $\lambda > 1$ and $C>0$, satisfies 
the inequality
\[
\frac{1}{\lambda}d(g,h)-C \leq d' (\phi(g),\phi(h)) 
\leq \lambda d(g,h) + C\,, \quad\text{ for all } g,h\in\Ga\,. 
\]
Throughout the paper, when a finite generating set is given for a group, 
the group shall always be equipped with the word metric, even when no 
metric is specified. We note that the word metric is quasi-isometrically 
invariant under different choices of finite generating set.

The aim of this paper is to prove that each member of a large class of 
right-angled Artin groups, which we shall call {\it planar type} 
right-angled Artin groups, embeds {\it quasi-isometrically} in some 
pure braid group $PB_m$, as well as in the group $\DiffD$ of 
area-preserving diffeomorphisms of the unit disk, equipped with the 
so-called ``hydrodynamical'' or $L^2$-norm metric. 
Quite independently one may observe 
that many interesting groups may be embedded quasi-isometrically in 
right-angled Artin groups. As a corollary, then, we obtain that all 
surface groups, with the exception of the three simplest non-orientable 
surfaces, as well as at least one hyperbolic 3-manifold group, embed 
quasi-isometrically in $PB_m$, for some $m$, and in $\DiffD$.

In \cite{CW} the authors showed that each right-angled Artin group of 
planar type embeds in a pure braid group $PB_m$, for some $m$ depending 
on the defining graph $\Delta$. However, it was not clear whether this 
embedding is quasi-isometric. In the present paper we modify the 
construction in order to give quasi-isometric
embeddings. Our techniques also yield quasi-isometric embeddings of 
arbitrary right-angled Artin groups into closed surface
mapping class groups. (See Theorem \ref{MainThm} and Corollary \ref{MainCor}).

The initial motivation for the present paper, however, came from the 
work \cite{BG} of Benaim and Gambaudo, who introduce the hydrodynamical 
metric $d_{\rm hydr}$ on the group $\DiffD$ of volume preserving 
diffeomorphisms of the closed disk $D^2$ (see Section \ref{S:PBtoDiff} 
for details). Observing that the metric $d_{\rm hydr}$ is 
unbounded they proposed a study of the large-scale properties of $\DiffD$
with respect to this metric. Their main result in this direction states 
that, for any $n$, the free abelian and free groups, $\Z^n$ and $F_n$, 
embed quasi-isometrically in $\DiffD$. Adapting their techniques to the 
case of right-angled Artin groups we are able to show that all of the 
examples referred to above (planar type right-angled Artin groups, 
surface groups, and other hyperbolic group examples) may also be 
embedded quasi-isometrically in $\DiffD$ with respect to $d_{\rm hydr}$.

It should also be stressed that we do not actually
construct a quasi-isometric embedding of $PB_m$ in $\DiffD$. In fact,
it is unknown whether $PB_m$ can be embedded as a subgroup of $\DiffD$
or $\Homeo$, and it is rather unlikely that this is possible in any 
natural way. For instance, it follows from the work of Morita \cite{Morita}
that there is no group-theoretic section 
of the natural homomorphism $\P_m\to PB_m$ where $\P_m<\DiffD$ denotes 
a subgroup of the area-preserving diffeomorphisms which are fixed on a
given set of $m$ disjoint closed disks in the interior of $D^2$.

The plan of this paper is as follows. 
In section \ref{S:pRAAGinPB} we define planar right-angled Artin groups,
and prove that they embed quasi-isometrically in a pure braid group
$PB_m$ for large enough $m$. This proof is the heart of the paper.
In section \ref{S:PBtoDiff} we deduce a quasi-isometric embedding of
planar right-angled Artin groups in $\DiffD$. This section will not come
as a surprise to any reader of \cite{BG}. Finally, in section 
\ref{S:HYPinRAAG} we build on the work in \cite{CW} in order find
many interesting quasi-isometrically embedded subgroups of $PB_m$
and $\DiffD$.


\section{Planar right-angled Artin groups in pure braid groups}
\label{S:pRAAGinPB} 


In this section we define the notion of ``planarity'' for right-angled 
Artin groups and show that a right-angled Artin group of planar type 
$G(\Delta)$ may be embedded quasi-isometrically in an $m$-strand pure braid 
group $PB_m$ (where $m$ depends on the defining graph $\Delta$) -- see 
Corollary \ref{MainCor}~(2). More generally, our techniques show that any 
right-angled Artin group (not necessarily of planar type) may 
be embedded quasi-isometrically in the mapping class group $\Mod(S)$ of 
an orientable punctured surface, surface with boundary or closed surface 
(with implicit restrictions on the genera in each case) -- see Theorem 
\ref{MainThm} and Corollary \ref{MainCor}~(1). In achieving these results 
we do not overly concern ourselves with the problem of minimising 
the number $m$ of stands in the target pure braid group or the genera, 
or number of punctures or boundary components, of the surface $S$ when 
the target is the mapping class group $\Mod(S)$. It would nevertheless 
be interesting to further understand for which $m$ the group $PB_m$, and 
for which surfaces $S$ the group $\Mod(S)$, admits a (quasi-isometric) 
embedding of a given right-angled Artin group. We first fix some notation.

\begin{notation}{\bf($\Mod(S)$, $\PMod(S)$) }
Let $S$ denote a (not necessarily connected)
\emph{finitely punctured compact orientable surface}:
$S$ is homeomorphic to $S_0\setminus P$ where $S_0$ is a compact 
orientable surface with boundary $\partial S_0$, and $P$ is a finite set 
of points in $S_0\setminus\partial S_0$. We note that the boundary of $S$ 
is just that of $S_0$, namely $\partial S=\partial S_0$.

We shall write simply $\Mod(S)$ for the orientable mapping class group of 
$S$ relative to the boundary. That is,  
the diffeomorphisms of $S$ are supposed to fix the boundary pointwise, 
and they are equivalent if they are related by a 
diffeotopy that is the identity on the boundary. We shall also write 
$\PMod(S)$ for the subgroup of $\Mod(S)$ generated by Dehn twists. 
This is simply the finite index normal subgroup of $\Mod(S)$ whose 
elements leave invariant the components of $S$ and do not permute the 
elements of the puncture set $P$ \cite{Dehn}.

Note that if $S$ is a disjoint union of surfaces 
$S_1$ and $S_2$ then we simply have $\Mod(S)\cong\Mod(S_1)\times \Mod(S_2)$.
\end{notation}  

\subsection{Planarity, circle diagrams, and the basic representation}
\label{subsect:definitions} 

Any finite collection $\Cal C=\{C_1,C_2,\ldots,C_n\}$ 
of connected compact subsets of $\R^2$ (or $D^2$ or $S^2$) determines a 
\emph{non-incidence} graph $\Delta_{\Cal C}$, defined to be the 
simplicial graph with vertex set $\Cal C$ and edges $\{C_i,C_j\}$ 
whenever $C_i\cap C_j=\emptyset$.

Note that any collection of plane diffeomorphisms $\phi_1,\ldots,\phi_n$ 
with $\phi_i$ trivial outside a regular neighbourhood of
the set $C_i$, for each $i=1,\ldots,n$, will generate a homomorphic image 
of the right-angled Artin group $G(\Delta_{\Cal C})$ 
associated to the graph $\Delta_{\Cal C}$.

\begin{defn}[\bf(Circle diagram; planar type.)]
We say that the right angled Artin group $G(\Delta)$ is of 
\emph{planar type} if $\Delta$ is isomorphic to the non-incidence graph 
$\Delta_{\Cal C}$ where $\Cal C$ denotes a finite collection of smooth 
simple closed curves in general position in the interior of the disk $D^2$. 
We call $\Cal C$ a \emph{(planar) circle diagram} for $G(\Delta)$.  

More generally, we may define a \emph{circle diagram} for $G(\Delta)$ to 
be a collection $\Cal C$ of simple closed curves in some orientable 
surface whose non-incidence graph $\Delta_{\Cal C}$ is isomorphic to 
$\Delta$. Note that one may easily find a circle diagram, in this larger 
sense, for an arbitrary right-angled Artin group (c.f. \cite{CW}).
\end{defn}

If $\Delta$ is a simplicial graph then we define its \emph{complementary} 
(or \emph{opposite}) graph $\Delta^{\text{op}}$ to be the simplicial
graph with the same vertex set as $\Delta$ and which has an edge between 
two vertices if and only if $\Delta$ does not.
We recall from \cite{CW} that if the complementary defining 
graph $\Delta^{\text{op}}$ is planar then $G(\Delta)$ is of planar type. 
We also recall the idea of the proof. 
An embedding of the graph $\Delta^{\text{op}}$ in the plane $\R^2$, 
with vertex set $\Cal V\subset \R^2$ gives rise to a collection 
of simple closed curves $\Cal C =\{C_v : v\in \Cal V\}$ in $\R^2$,
where $C_v$ is defined as the boundary of a regular neighbourhood
of the union of $v$ and one-half of each edge adjacent to $v$.

We stress that the planarity of $\Delta^{\rm op}$ is a sufficient, but
by no means necessary condition for $G(\Delta)$ to be of planar type. 
This point is nicely illustrated by the example described in 
Section \ref{subsect:3mfld}. Figure \ref{crazycurves} shows a planar 
circle diagram for a graph (the 1-skeleton of an icosahedron) whose 
complementary graph is non-planar.

\begin{figure}[htb] 
\begin{center}
\input{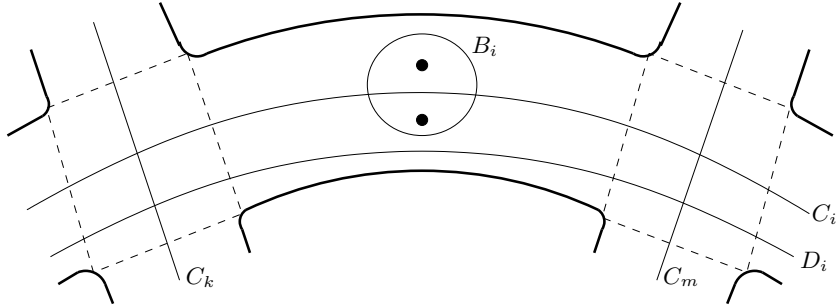}
\end{center}
\caption{A segment of the $i$th annulus $A_i$, and its intersection
with two other annuli, $A_k$ and $A_m$. The curves $B_i,C_i,D_i$ are 
indicated with solid lines, and the boundaries of the squares $A_i\cap A_k$ 
and $A_i\cap A_m$ with dashed lines.}\label{F:cuttingarcs}
\end{figure}

\begin{defn}[\bf(Surface associated to a circle diagram.)]\label{DefnSC}
Given a (smooth) circle diagram $\Cal C=\{C_1,\ldots,C_n\}$ 
(in an arbitrary orientable surface) we 
define a compact surface with boundary $S_{\Cal C}$ associated 
to $\Cal C$, as follows.

Let $S'$ denote a regular closed neighbourhood of the circle diagram. 
Thus $S'$ is a union of annuli $A_i$ (with $A_i$ a regular neighbourhood 
of $C_i$). Moreover, each intersection point of the curves $C_i$ and $C_j$ 
gives rise to one square in the surface $S'$, which is just one path 
component of $A_i\cap A_j$. The whole surface $S'$ is a compact
orientable surface with boundary.

In each annulus $A_i$ we introduce a pair of distinguished points 
which do not lie in the intersection with another annulus. The two 
points must be on opposite sides of the curve $C_i$, as indicated in 
figure \ref{F:cuttingarcs}. We denote $P$ the union of all these 
distinguished points. Finally, we define $S_{\Cal C}$ to be the surface
\[
S_{\Cal C} := S'\setminus N(P)\,,
\]
where $N(P)$ denotes a regular open neighbourhood of the finite set $P$.  
\end{defn}

We remark that in the preceding construction two annuli $A_i$ and $A_j$
are disjoint if and only if the corresponding generators $a_i$ and $a_j$
of the right-angled Artin group $G(\Delta_\mathcal{C})$ commute.

\begin{defn}[\bf(Basic representation 
$G(\Delta_{\Cal C})\to \Mod(S_{\Cal C})$.)]
To a circle diagram $\mathcal{C}=\{C_1,\ldots,C_n\}$ in an orientable
surface we can associate a representation 
$G(\Delta_{\Cal C})\to \Mod(S_{\Cal C})$ (from the right-angled Artin 
group whose defining graph is the non-adjacency graph of $\mathcal{C}$
to the mapping class group of the surface $S_\mathcal{C}$) as follows.
In each annulus $A_i$ of $S$ we draw smooth simple closed curves 
$B_i, C_i$ and $D_i$ as indicated in figure \ref{F:cuttingarcs}, and 
define the following diffeomorphism for each $i=1,\ldots,n$:
\[
f_i=\tau_{B_i}\circ\tau_{D_i}^2\circ\tau_{C_i}^{-2}\circ
\tau_{B_i}\in \Diff(S_{\Cal C},\partial S_{\Cal C})\,,
\]
where $\tau_C$ denotes a smooth Dehn twist along a curve $C$.
(We remark that the diffeomorphism $\tau_{D_i}^{-2}\circ\tau_{C_i}^2$
may be thought of as induced by the pure braid of the set $P$
of marked points on $S'$ which is given by moving the puncture enclosed by
the curves $C_i$ and $D_i$ twice in an anticlockwise sense around the 
annulus.)

This clearly defines a homomorphism 
$f\co G(\Delta_{\Cal C})\to\Diff(S_{\Cal C},\partial S_{\Cal C})$ by 
setting $f(a_i)=f_i$.

Whenever $S_{\Cal C}$ is viewed as a subsurface of any other (not 
necessarily compact) surface $\what S$ the homomorphism $f$ extends 
naturally to a homomorphism 
\[
\what f\co G(\Delta_{\Cal C})\to \Diff(\what S, \partial \what S)\,,
\]
where every element of the image acts by the identity on 
$\what S\setminus S$.

We shall denote $\varphi\co G(\Delta_{\Cal C})\to \Mod(S_{\Cal C})$ and 
$\what\varphi\co G(\Delta_{\Cal C})\to \Mod(\what S)$ the homomorphisms
induced by $f$ and $\what f$ respectively. Clearly $\what\varphi$ is 
obtained from $\varphi$ by composing with the map 
$\Mod(S_{\Cal C})\to \Mod(\what S)$ induced by the inclusion.
\end{defn}

We remark that in the above construction, the homeomorphism $f_i$,
when restricted to the annulus $A_i$, is pseudo-Anosov. This idea,
which is essential to our proof that the embedding in quasi-isometric,
is inspired by \cite{DynnWie}.

Note also that throughout the whole of the above discussion we did not need 
to suppose that any of the surfaces are necessarily connected. If $S$ is 
a disjoint union of surfaces $S_1$ and $S_2$ then we understand 
$\Mod(S)\cong\Mod(S_1)\times \Mod(S_2)$.

\subsection{Quasi-isometric embeddings}

Our next aim is to prove that the basic representation 
$G(\Delta_{\Cal C})\to \Mod(S_{\Cal C})$ is injective and quasi-isometric.
In fact, we are going to prove something slightly stronger, which is the
main technical result of this section:

\begin{thm}\label{MainThm}
Let $\Cal C$ be a circle diagram, $\Delta_{\Cal C}$ the associated 
non-incidence graph, and $S_{\Cal C}$ the surface of Definition 
\ref{DefnSC}. Suppose that $S_{\Cal C}$ embeds as a subsurface of an 
orientable finitely punctured compact surface $\what S$ and that the 
embedding $S_{\Cal C}\hookrightarrow \what S$ is $\pi_1$-injective on 
each component of $S_{\Cal C}$. Then  the homomorphism 
\[
\what\varphi\co G(\Delta_{\Cal C})\to \PMod(\what S)
\]
is a quasi-isometric embedding of groups (in particular, an injective 
homomorphism).
\end{thm}

Before passing on to the proof of this theorem, we mention some easy 
consequences. Firstly, since it is a finite index subgroup, the inclusion 
of $\PMod(\what S)$ into $\Mod(\what S)$ is a quasi-isometric embedding.
It follows from the statement of Theorem \ref{MainThm} that we also obtain
a quasi-isometric embedding of $G(\Delta_{\Cal C})$ into the slightly 
larger group $\Mod(\what S)$.

\begin{cor}\label{MainCor}
Let $G=G(\Delta)$ denote a right-angled Artin group.
\begin{enumerate}
\item The group $G$ embeds quasi-isometrically in $\Mod(S)$ for some 
connected closed orientable surface $S$ (of genus depending on $\Delta$). 
\item If $G$ is of planar type then it embeds quasi-isometrically in the 
pure braid group $PB_m$ (for a sufficiently large $m$ depending on 
$\Delta$).
\end{enumerate}
\end{cor}

\begin{proof} 
For an arbitrary right-angled Artin group $G(\Delta)$ we may always find a 
circle diagram $\Cal C$ on some orientable surface such that
$\Delta=\Delta_{\Cal C}$. Suppose that $S_{\Cal C}$ has $b$ boundary 
components. Then we define a $\pi_1$-injective inclusion  of $S_{\Cal C}$ 
into a closed connected surface $\what S$ by gluing $S_{\Cal C}$ along 
its boundary to an orientable genus zero surface with $b$ boundary 
components. Statement (1) now follows by applying Theorem \ref{MainThm}.

When $G(\Delta)$ is of planar type we may find a planar circle diagram 
$\Cal C$ with $\Delta=\Delta_{\Cal C}$. The corresponding surface 
$S_{\Cal C}$ is then also planar, and may be viewed as a subsurface of 
$D^2$. Removing a single point from each disk component of 
$D^2\setminus S_{\Cal C}$ yields an $m$-punctured closed disk 
$\what S\cong D^2\setminus\{\, m \text{ points}\,\}$ (such that 
$\partial \what S\cong S^1$). 
The inclusion $S_{\Cal C}\hookrightarrow \what S$ is, by construction, 
$\pi_1$-injective on each connected component of $S_{\Cal C}$.
We also recall the fact that $\Mod(\what S)$ is naturally isomorphic to 
the $m$-string braid group $B_m$, for some $m$, and 
$\PMod(\what S)\cong PB_m$, the pure braid group. 
Statement (2) now follows from Theorem \ref{MainThm}.
\end{proof}

\subsection{Proof of Theorem \ref{MainThm}}

We suppose throughout that $\Delta=\Delta_{\Cal C}$ where 
$\Cal C=\{C_1,\ldots,C_n\}$ is a smooth circle diagram in some orientable 
surface, and $S=S_{\Cal C}$ the compact orientable surface of Definition 
\ref{DefnSC}. For simplicity (and without any loss of generality) we shall 
view $\Cal C$ as being a circle diagram in the surface $S$.
We assume that an inclusion $S\to \what S$ is given, and that  
the maps $f,\what f,\varphi,\what\varphi$ are as described in the 
preceding definitions (Subsection \ref{subsect:definitions}).

\begin{defn}{\bf (Curve diagrams and intersection numbers)}
\label{D:CurveDiag} By a \emph{(smooth) curve diagram} on a surface we 
shall mean the diffeotopy class of the
union of a collection of (not necessarily disjoint) smooth simple 
closed curves on the surface, no two of which are homotopy equivalent. 

If $D,D'$ are two curve diagrams then we write $|D\cap D'|$ for the 
minimal intersection number between representatives of the 
diffeotopy classes $D$ and $D'$ which are chosen to be in 
transverse position with respect to one another.
\end{defn}

We observe that the group $\Mod(S)$ acts naturally (on the left) on the
set of curve diagrams.

We shall first investigate the action of $G(\Delta)$ on the surface 
$S=S_{\Cal C}$ via the homomorphism $\varphi\co G(\Delta)\to\Mod(S)$. 
We shall let $a_1,\ldots,a_n$ denote the standard generators of 
$G(\Delta)$. Recall that $S$ is the union of annuli $A_i$, $i=1,\ldots,n$,
with a pair of open disks removed from each annulus -- see Subsection 
\ref{subsect:definitions}.

We define an $(i,k)$-{\it intersection square} to be one connected 
component of $A_i\cap A_k$, the intersection of the $i$th and the $k$th 
annulus.

When drawing pictures of curve diagrams we shall always arrange for the
curves to be {\it reduced} with respect to the four boundary lines of 
each intersection square -- this means that there are no bigons enclosed
between the curves and the boundary of the intersection square. This
is equivalent to asking that the number of intersections between the
curves and the intersection squares be minimal in the diffeotopy class
of the curve diagram.


\begin{figure}[htb] 
\begin{center}
\input{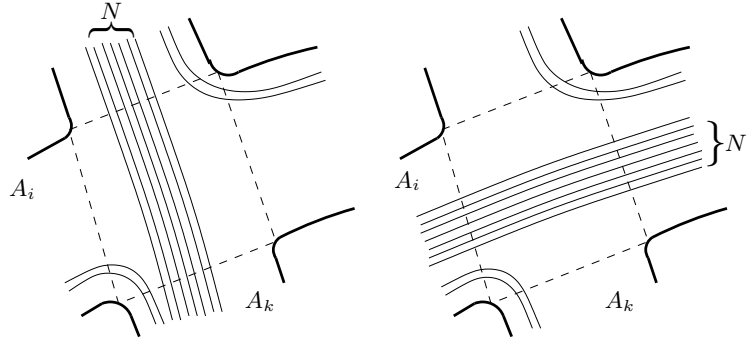}
\end{center}
\caption{A curve diagram $D$ traverses an $(i,k)$-square $N$ times in the 
$A_k$-direction, and $N$ times in the $A_i$ direction.}
\label{F:intersectsquare}
\end{figure}


If $Q$ is an $(i,k)$-intersection square, then we say that a curve diagram 
$D$ traverses the square $N$ times {\it in the $A_i$-direction} 
(resp.\ \emph{in the $A_k$-direction}) if, in a version of $D$ which
is reduced with respect to the intersection squares, the following
condition is satisfied: among the connected components of $Q\cap D$ 
there are precisely $N$ which connect opposite sides of the square 
without leaving the interior of $A_i$ (resp. $A_k$).
See figure \ref{F:intersectsquare}.

\begin{notation}
Let $D$ be a curve diagram in $S$. For each $i=1,\ldots,n$, write $c_i(D)$ 
for the largest number $N$ such that $D$ traverses every $(i,k)$-square on 
the annulus $A_i$ at least $N$ times in one or other direction. 
We say that $D$ is \emph{transverse} to $A_i$ if there exists an
$(i,k)$-square (where $k\in \{1,\ldots,n\}$ is such that $a_i$ 
and $a_k$ do not commute) which is traversed by $D$ 
at least $c_i(D)$ times in the $A_k$-direction.
\end{notation}

\begin{lemma}\label{L:ComplIncrease}
Let $D$ be a curve diagram in $S$, $a_i$ some generator of $G(\Delta)$, and
let $D'=\varphi(a_i)^p(D)$. Suppose, furthermore, that the diagram $D$ is 
transverse to the annulus $A_i$. Then $c_i(D')\geqslant 2^{|p|}c_i(D)$. 
Moreover, $D'$ is transverse to every annulus $A_k$ which intersects $A_i$ 
(and not transverse to $A_i$).
\end{lemma}

\begin{figure}[htb] 
\begin{center}
\input{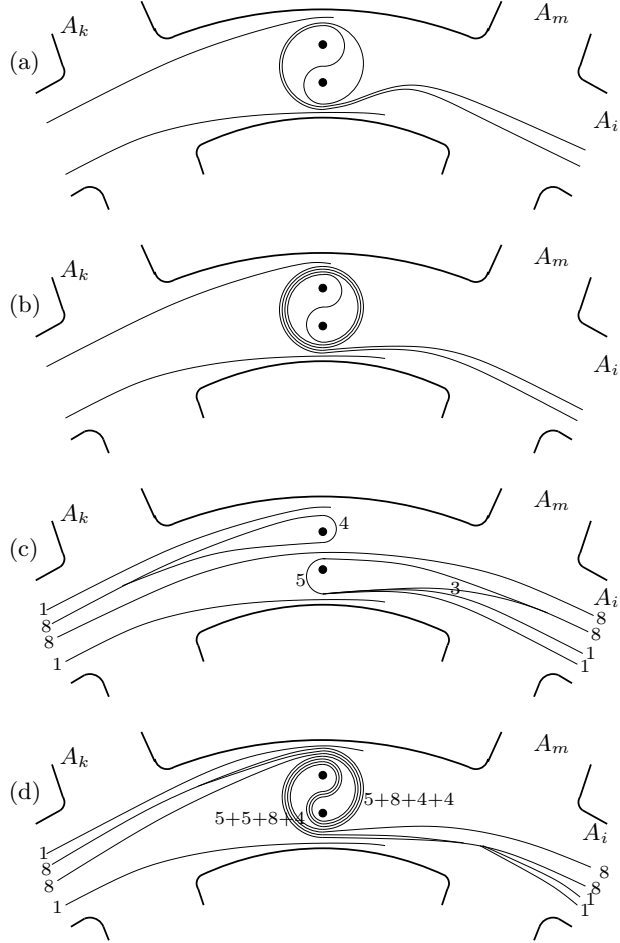}
\end{center}
\caption{Acting by $\what f(a_i)$ twice in a row. (a) After the first 
action by $\what f(a_i)$. \ (b) After the action by $\tau_{B_i}$. \ 
(c) After the action by $\tau_{C_i}^2\tau_{D_i}^{-2}$. \ (d) After the
action by $\tau_{B_i}$.}\label{F:case2}
\end{figure}

\begin{proof}
Note that, since $D$ is transverse to $A_i$ there is some 
$(i,k)$-intersection square where $D$ traverses at least $c_i(D)$ 
times in the $A_k$ direction. We claim that the image of each of these 
crossing arcs under the action of $\varphi(a_i)$ will traverse each 
intersection square on the annulus $A_i$ at least $2^p$ times in the 
$A_i$ direction (once $D'$ is placed in minimal position with respect 
to the boundary arcs of the intersection squares). 
This clearly establishes the Lemma. 
(To see that $D'$ is not transverse to $A_i$ observe that if it were then, 
since $D=\varphi^{-p}(a_i)(D')$, a further application of the Lemma would 
imply that $c_i(D)\geqslant 2^{|p|}c_i(D')\geqslant 2^{|p|}c_i(D)$, a 
contradiction).

Let us first consider the special case $p=1$, that is, the action
by $\varphi(a_i)$ on a single crossing arc $\alpha$ (the case $p=-1$ is 
of course similar). We first recall that $\varphi(a_i)$ is defined by the 
diffeomorphism
\[
f_i=\tau_{B_i}\circ\tau_{D_i}^2\circ\tau_{C_i}^{-2}\circ
\tau_{B_i}\in \Diff(S,\partial S)\,.
\]
Figure \ref{F:case1} illustrates the different stages 
of the action of this element on $\alpha$. 
At the final stage the image of $\alpha$ traverses each 
intersection square at least twice.


\begin{figure}[ht] 
\begin{center}
\input{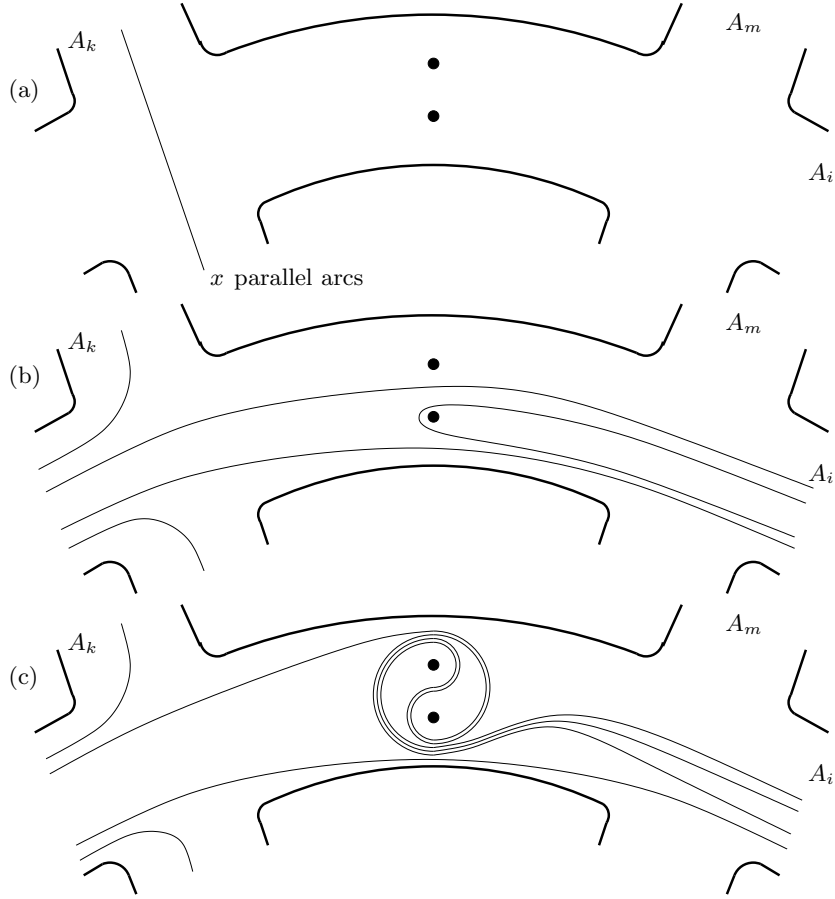}
\end{center}
\caption{If there are $x$ arcs traversing a square of $A_i\cap A_k$
in the $A_k$-direction, then after the action of $\what f(a_i)$ there
are $2x$ arcs traversing any square of intersection involving $A_i$
in the $A_i$-direction. (Part (b) shows the result of the action by 
$\tau_{B_i}\circ \tau_{D_i}^2\circ \tau_{C_i}^{-2}$.)
}
\label{F:case1}
\end{figure}


Now we consider the more general case $p\in \N$ (and again, the case
$-p\in \N$ is similar). We have just seen that the image of $\alpha$
under $\varphi(a_i)$ contains the segment 
$\alpha'$ shown in Figure \ref{F:case2}(a). Now Figure
\ref{F:case2} shows that acting once more on this curve $\alpha'$ by
$\varphi(a_i)$ yields a curve diagram containing at least twice (in
fact: five times) as many parallel copies of the same segment $\alpha'$ 
appearing as subsegments. Inductively, the curve diagram of
$\varphi(a_i)^p(\alpha)$ contains at least~$2^{p-1}$
copies of this arc 
$\alpha'$, and hence traverses each intersection square in $A_i$ at least 
$2^p$ times in the $A_i$ direction.
\end{proof}


We shall say that a word in the letters $a_1^{\pm 1},\ldots,a_n^{\pm 1}$ 
is {\it reduced} if there is no shorter words in those letters
representing the same element of the right-angled Artin group $G(\Delta)$.
It is well-known that any two reduced words differ by finite a sequence 
of ``shuffles'': exchanges of adjacent letters $a_i^\pm$, $a_j^\pm$ where 
$i,j$ span an edge in $\Delta$ (so that $a_i,a_j$ commute). This seems to 
be first due to Baudisch \cite{Baud}.
A much more recent proof can be found in \cite{CW} (Proposition 9(i)).
It follows from this result that, for each $i=1,..,n$, the number $\ell_i$ 
of occurrences of the letters $a_i$ and $a_i^{-1}$ in any reduced
representative is an invariant of the group element. 

\begin{lemma}\label{cjlemma}
Let $w\in G(\Delta)$, $\phi=\varphi(w)\in \Mod(S)$, and $D$ a curve diagram 
in $S$. Suppose that $D$ is transverse to $A_i$ for every $i$ for which $w$ 
can be written $w=v a_i^{\pm 1}$ in reduced form. Then, for all $j$, 
\[
c_j(\phi(D))\geqslant 2^{\ell_j(w)}c_j(D)\,,
\]
where $\ell_j(w)$ denotes the number of occurrences of the letters 
$a_j,a_j^{-1}$ in some, or any, reduced word representing $w$.
\end{lemma}

\begin{proof}
Suppose $w=v a_i^p$ in reduced form, where $|p|$ is as large as 
possible. Then if $v=u a_j$ in reduced form for some $j$, we  have 
$i\neq j$. Let $D'=\varphi(a_i)^p(D)$.  Since $D$ is transverse to $A_i$
it follows, by Lemma \ref{L:ComplIncrease}, that 
$c_i(D')\geqslant 2^{|p|}c_i(D)$ and, moreover, that $D'$ is transverse to 
every annulus $A_k$ which intersects $A_i$. 

Recall that $\varphi(a_i)$ is represented by a diffeomorphism $f_i$
which is the identity outside $A_i$. It follows that applying $f_i$ to 
any diagram does not increase the number of crossings (in either direction) 
at any $(j,k)$-square for $j\neq i$ (and any $k$ not adjacent to $j$).
Since $f_i^{-1}$ cannot increase these numbers either (for  the same 
reason), we conclude that the diffeomorphisms $f_i,f_i^{-1}$ can neither 
increase nor decrease the number of such crossings. Thus $c_j(D')=c_j(D)$ 
for $j\neq i$. Moreover, if $D$ is transverse to $A_j$ and $A_j$ is 
disjoint from $A_i$ then $D'$ is also transverse to $A_j$.

We now claim that the transversality hypothesis (on $D$ and $w$) applies 
once again to the element $v\in G(\Delta)$ and the diagram $D'$. For, if 
$v=u a_j$ in reduced form, then $j\neq i$ and either 
\begin{enumerate}
\item   $i,j$ are non-adjacent in $\Delta$: here $A_j$ intersects $A_i$ 
and so $D'$ is transverse to $A_j$, by the above application of Lemma 
\ref{L:ComplIncrease}; or
\item   $i,j$ span an edge of $\Delta$: here $A_j$ and $A_i$ are disjoint. 
In this case we have that $w=u a_j a_i=u a_i a_j$ in reduced form, so that 
$D$ is transverse to $A_j$ (by hypothesis). As observed above, in this 
case $D'$ remains also transverse to $A_j$.
\end{enumerate}

The Lemma now follows by applying the statement inductively to the pair 
$v$, $D'$.
\end{proof}

Recall that the curves of the original circle diagram 
$\Cal C=\{ C_1,\ldots,C_n\}$ appear as simple closed (but not disjoint)
curves in $S$, and since no two are homotopic, their union 
$E=\bigcup\Cal C$ is a curve diagram in $S$ -- indeed, we can think of $E$
as the trivial curve diagram in $S$. We note that $E$ is 
transverse to every annulus, and that $c_j(E)=1$ for each $j=1,\ldots,n$.
By Lemma \ref{cjlemma}, this latter fact characterises $E$ amongst 
all its translates by elements of the group $G(\Delta)$ (acting via 
$\varphi$). As a consequence, $\varphi\co G(\Delta)\to\Mod(S)$ is an 
injective homomorphism.

We now return to studying the homomorphism 
$\what\varphi\co G(\Delta)\to\PMod(\what S)$.
Our first observation is the following:

\begin{prop}\label{faithful}
The homomorphism $\what\varphi\co G(\Delta)\to\PMod(\what S)$ is faithful.
\end{prop}

\begin{proof}
Since the inclusion $S\hookrightarrow \what S$ is $\pi_1$-injective on 
each component, the numbers $c_i(D)$ are invariants of the homotopy (or 
isotopy) class of the diagram $D$ in $\what S$, just as in $S$. Now let 
$D=\what\varphi(a)(E)$. By Lemma \ref{cjlemma}, and the fact that $E$ 
is transverse to every annulus, $c_j(D)\geqslant 2$ for some $j$ unless 
$a=1$. This proves injectivity of the homomorphism $\what\varphi$. 
\end{proof} 

We note that injectivity of $\what\varphi$ is also a 
consequence of the statement that we are about to prove, namely that
$\what\varphi$ satisfies a quasi-isometric inequality. 
This is simply because $G(\Delta)$ is a torsion free group while, in 
general, any quasi-isometric homomorphism must have finite kernel.
\medskip

Let $\what{\Cal C}$ denote a finite collection of essential simple closed 
curves in $\what S$ which include the set $\Cal C =\{ C_i:i=1,\ldots,n\}$ 
and such that the Dehn twists along these curves generate $\PMod(\what S)$.
Observe that since the inclusion $S\to\what S$ is $\pi_1$-injective, no 
two curves of $\Cal C$ are homotopic in $\what S$ (note also that no curve 
$C_i$ can ever be parallel to a boundary component of $S$). We may 
therefore choose $\what{\Cal C}$ to be minimal in the sense that no two 
curves are homotopic. It follows that their union 
$\what E=\bigcup \what{\Cal C}$ is a curve diagram in $\what S$ 
(c.f.~Definition \ref{D:CurveDiag}). We note that $E\subset \what E$ is a 
subdiagram. We can think of $\what E$ as the extended trivial curve 
diagram in $\what S$.

We shall henceforth also fix the set of Dehn twists 
$\{ \tau_C \co C\in\what{\Cal C}\}$ as our choice of generating
set for $\PMod(\what S)$, and we shall write $d$ for the word metric in 
$\PMod(\what S)$ with respect to these generators. 

\begin{defn}{\bf (Complexity of an element of $\PMod(\what S)$)}
Suppose that the generating curves $\what{\Cal C}$ are chosen as above, 
so that $\what E$ is a curve diagram in $\what S$. If $\phi$ is an element
of the group $\PMod(\what S)$ then $\phi(\what E)$ denotes the curve diagram 
which is represented by the image of $\what E$ under any diffeomorphism 
representing $\phi$. We define the {\it complexity} of a curve diagram $D$ 
in $\what S$ to be 
\[
\comp(D)=\log_2(|D\cap \what E|) - \log_2(|\what E\cap \what E|)\,.
\] 
The complexity of an element $\phi\in\Mod(\what S)$ is defined by 
\[\comp(\phi)=\comp(\phi(\what E))\,.\]  
Note that the definition is normalized so that $\comp(\what E)=0$. 
\end{defn}

\paragraph{\bf Proof of Theorem \ref{MainThm}\ }
Theorem \ref{MainThm} is now a consequence of the faithfulness 
of $\what\varphi$ established in Proposition \ref{faithful}, and the 
following two propositions. Together, Propositions \ref{P:firstkey} 
and \ref{P:secondkey} imply that for an element $a\in G(\Delta)$ of 
wordlength $\ell(a)$ we have
\[
\ell(a)\leqslant K_1\cdot\comp(\what \varphi(a)) +K_0 
\leqslant K_1K_2\cdot d(\what \varphi(a),1) +K_0,
\]
from which it follows that the homomorphism $\what\varphi$ is a 
quasi-isometric embedding.

\begin{prop}\label{P:firstkey}
Let $\ell(a)$ denote the length of a shortest representative word of an 
element $a\in G(\Delta)$. Then the complexity of $\what \varphi(a)$ 
grows {\it at least} linearly with $\ell(a)$.  In other words
\[
\ell(a)\leqslant K_1\cdot\comp(\what \varphi(a)) + K_0\,.
\]
where $K_0$ and $K_1$ are some positive constants (e.g. $K_1$ equal to the 
number of generators in $G(\Delta)$, and $K_0=K_1\cdot \log_2(|\what E\cap \what E|)$ suffice).
\end{prop}

\begin{proof}
Recall that the subdiagram $E\subset \what E$ is a diagram in $S$ which is
transverse to every annulus, and that $c_j(E)=1$ for each $j$.
Let $D=\varphi(a)(E)$. Then it follows, by Lemma \ref{cjlemma}, 
that $c_j(D)\geqslant 2^{\ell_j(a)}$, for each $j=1,\ldots,n$.

Suppose that $C,C'$ are essential simple closed curves in $S$. Then we 
observe that, since the inclusion $S\to\what S$ is $\pi_1$-injective, 
the intersection number $|C\cap C'|$ will be the same whether measured 
in $S$ or in $\what S$. (This is because, if $C$ and $C'$ cobound a bigon 
in $\what S$ then they must also cobound a bigon in $S$.) It follows by 
this reasoning that, for each $j=1,\ldots,n$,
\[
c_j(D)\leqslant |D\cap E|\leqslant |\what D\cap \what E|\, ,
\]
where $\what D:=\what\varphi(a)(\what E)$ and $D:=\what\varphi(a)(E)$. Thus
\[
2^{\ell(a)/n}\leqslant \max\{2^{\ell_j(a)}:j=1,\ldots,n\}\leqslant
 \max\{c_j(D):j=1,\ldots,n\}\leqslant|\what D\cap \what E|\,,
\]
and the proposition follows with $K_1=n$ and 
$K_0=K_1\cdot\log_2(|\what E\cap \what E|)$.
\end{proof}

\begin{prop}\label{P:secondkey}
The complexity of the curve diagram of an element $\phi$ of $\PMod(\what S)$
grows at most linearly with the distance of $\phi$ from the neutral
element in the Cayley graph of $\PMod(\what S)$ -- that is, we have 
\[
\comp(\phi)\leqslant K_2\cdot d(\phi,1),
\]
where $K_2$ is a positive constant (equal to the base 2 logarithm of the
number of curves in $\what{\Cal C}$).
\end{prop}

\begin{proof}
We take the Dehn twists along the curves in $\what{\Cal C}$ as our finite 
generating set for $\PMod(\what S)$. If $\tau_C$ is a given generator 
(the Dehn twist along the curve $C\in\what{\Cal C}$) and $\what D$ a curve 
diagram of known complexity, then we can easily estimate the complexity of 
$\tau_C(\what D)$. Namely, for each point of intersection of $\what D$ with 
$C$, application of $\tau_C$ may introduce at most $r$ new points of
intersection with curves in $\what E$, where 
\[
r= \#\{ \text{ curves of } \what{\Cal C} \text{ which intersect } C\,\} 
\leqslant N-1\,, 
\]
where $N=|\what{\Cal C}|$. 
Since $|C\cap \what D|\leqslant |\what E\cap \what D|$, we then have 
\[ 
|\what E\cap \tau_C(\what D)| \leqslant 
|\what E\cap \what D| + (N-1)\cdot |C\cap \what D| \leqslant 
N\cdot |\what E\cap \what D|\,. 
\] 
Thus $\comp(\tau_C(\what D)) \leqslant \comp(\what D)+\log_2(N)$. 
By a straightforward induction on the wordlength of $\phi$, we then obtain 
\[ 
\comp(\phi) \leqslant K_2\cdot d(\phi,1)\,, 
\] 
where $K_2=\log_2(N)$.
\end{proof}


\section{From pure braid groups to $\DiffD$}\label{S:PBtoDiff} 


In the previous section we proved that there exists a 
homomorphic quasi-isometric embedding of any right-angled Artin 
group of planar type $G$ into the pure braid group of a 
certain number of points $P_1,\ldots,P_m$ in a disk $D^2$. 
The number of points needed depends, of course, on the group's 
circle diagram. We shall assume that the disk
is the disk with radius 1 and centre $(0,0)$ in the plane,
and that the distinguished points are 
$P_i=(\frac{i}{m-1}-\frac{1}{2},0)$ (with $i=0,\ldots,m-1$).

We also recall that the homomorphism constructed in the last
section factors through a certain subgroup of 
$\textsl{Diff}(D^2,\partial D^2)$, namely the group of
diffeomorphisms of the disk which fix pointwise the distinguished
points.

In the current section we shall point out that there is, in fact
an embedding in the group $\mathcal{P}_m$ of \emph{volume-preserving}
diffeomorphisms of the disk which, moreover, fix pointwise not only
the distinguished points and the boundary of the disk, but even disks 
of radius $r$ centered on each of the distinguished points, as well
as a regular neighbourhood of the boundary; here $r$ is a sufficiently 
small positive real number. This yields an embedding of $G$ in the group
$\mathcal{P}_m$, which is itself a subgroup of the group $\DiffD$
of volume-preserving diffeomorphisms of the disk $D^2$ which are the
identity on a neighbourhood of $\partial D^2$.

The aim of the current section is to prove that, if we equip 
$\DiffD$ with the \emph{hydrodynamical metric}, then this homomorphism
$G\to \DiffD$ is itself quasi-isometric. (This is stronger than saying
that the composition $G\to \mathcal{P}_m \to 
PB(D^2, P_1\cup\ldots\cup P_m)$ is quasi-isometric.)

We recall from \cite{BG} the definition of the hydrodynamical metric:
if $\{\phi_t\}_{t\in[0,1]}$ is a path in $\DiffD$, then the length of
this path is 
$$
\int_0^1 \sqrt{\int_{D^2} \left\|\frac{d\phi_t}{dt}(x) \right\|^2 dx }\ dt,
$$
where the symbol $\|.\|$ denotes the Euclidean norm of a tangent
vector to the disk. The hydrodynamical length $l_\textrm{hydr}(\phi)$ 
of an element $\phi$ of $\DiffD$ is then the infimum length of a path 
from the identity map to $\phi$. This defines a left-invariant
metric $d_\textrm{hydr}$ by setting
$d_\textrm{hydr}(\phi,\psi)=l_\textrm{hydr}(\phi^{-1}\psi)$.

Now we recall a technical result of Benaim and Gambaudo (Lemma 4 in 
\cite{BG}). There exists a constant $K>0$ and a function $C\co \R_+\to\R_+$ 
such that $\lim_{r\to 0} C(r) = 0$ with the following property.
Suppose that $\phi$ is an element of $\mathcal{P}_m$ (fixing disks
of radius $r$ around each of the distinguished points), and which
represents an element of the pure braid group whose shortest expression
as a product of Artin's standard generators $\sigma_1^{\pm 1},
\ldots,\sigma_{m-1}^{\pm 1}$ has length $l_{\mathrm{Artin}}$, then
$$
l_\textrm{hydr}(\phi) \geqslant \frac{1}{K} \cdot l_{\mathrm{Artin}}
\cdot (1-C(r)) \cdot (\text{area}\ D_0)^2.
$$
Here $D_0$ denotes a disk of radius $r$. 

This technical result implies immediately the following theorem:

\begin{thm}[Benaim-Gambaudo \cite{BG}]\label{BGthm} Suppose that $f$ 
is a homomorphism from a group $G$ to $\mathcal{P}_m$
for some choice of $m$ fixed disks around the puncture points, 
and suppose that the induced homomorphism $\varphi\co G\to PB_m$ is 
a quasi-isometric embedding, then so is $f$. 
\end{thm}

We are now ready to prove one of the main results of this paper.

\begin{thm}
If $G$ is a planar right-angled Artin group, then there exists 
a quasi-isometric embedding of $G$ into $\DiffD$.
\end{thm}

\begin{proof}
Corollary \ref{MainCor}(ii) gives a homomorphism 
$\what\varphi \co G(\Delta)\to PB_m$ which, by the construction given in
Section \ref{subsect:definitions}, factors through a homomorphism 
$\what f \co G(\Delta)\to\Diff(S, \partial S)$ where $S$ is a compact 
subsurface of the $m$-punctured disk. We may suppose, in fact, that $S$ is 
just the closed disk with a $m$ open disks removed, and is contained in 
the exterior of a suitably chosen collection of $m$ open disks of some 
constant radius $r>0$. We claim that the diffeomorphisms $\what f(a_i)$ 
which generate the image of $\what f$ may be chosen to be area preserving. 
It then follows that the image of $\what f$ lies in the subgroup 
$\P_m <\DiffD$ (defined by the collection of disks of radius $r$ just 
mentioned). Since the induced map $\what \varphi\co G(\Delta)\to PB_m$ is 
shown to be a quasi-isometric embedding (Corollary \ref{MainCor}(ii)), 
the result now follows by Theorem \ref{BGthm}.

To justify the claim, it suffices to observe that any Dehn twist about 
a smooth curve $C$ in $D^2$ may be realised by a volume preserving 
diffeomorphism with support in an arbitrarily small neighbourhood of 
the curve $C$. This follows from the following two observations:

\begin{itemize}
\item[(i)]
a smooth embedding $c\co S^1\to D^2$ can be extended, for 
a sufficiently small $\ep > 0$, to a smooth area preserving embedding
$S^1\times [-\ep,\ep]\to D^2$.

\item[(ii)] 
if we choose a smooth function $h\co [-\ep,\ep]\to [0,2\pi]$ 
such that $h(x)=0$ for $x<\frac{-\ep}{2}$, and $h(x)=2\pi$ for 
$x>\frac{\ep}{2}$, then 
\[
T_h\co S^1\times [-\ep,\ep]\to S^1\times [-\ep,\ep] \text{ \ \ such that \ } 
T_h(t,s)=(t+h(s),s)
\]
\end{itemize}
defines a smooth area preserving diffeomorphism which is the identity 
outside $S^1\times [\frac{-\ep}{2},\frac{\ep}{2}]$.
\end{proof}


\section{Quasi-isometrically embedded hyperbolic subgroups}
\label{S:HYPinRAAG} 

In this section we consider examples of quasi-isometric embeddings 
of groups in right-angled Artin groups of planar type.
\subsection{Hyperbolic surface groups.}

As proved in \cite{CW} every closed hyperbolic surface group with 
$\chi(S)\leqslant -2$ embeds quasi-isometrically in some right-angled 
Artin group of planar type. Therefore each of these surface
groups embeds quasi-isometrically both in $PB_m$ for some $m$, and 
in $\DiffD$. 

The case of an orientable surface $F$ of genus 2 is illustrated in 
Figure \ref{F:extension}. The figure describes a homomorphism 
$\pi_1(F)\to G(C_5)$ where $C_5$ denotes the 5-cycle graph with generators 
$a,b,c,d,e$. The homomorphism is determined by fixing a basepoint $\star$ 
and, for each loop $\gamma$ based at $\star$ in $F$, reading the sequence 
of crossings with sign that $\gamma$ makes with the transversely labelled
curves on the surface $F$.
We remark that this homomorphism actually
projects to an embedding into the corresponding right-angled Coxeter
group $W(C_5)$. (Note that $W(C_5)$ acts by isometries of the hyperbolic plane and
that the surface subgroup obtained here is finite index in the Coxeter group).

Embeddings of arbritary higher genus orientable surface groups may be 
obtained by simply restricting to finite index subgroups of $\pi_1(F)$. 
The treatment of the non-orientable case is similar but more complicated to 
describe (see \cite{CW}, Section 4).

The fact that these embeddings are all quasi-isometric is a consequence of 
the method used in \cite{CW} to prove injectivity: namely, each embedding 
is realised as the homomorphism induced on the fundamental groups by a 
locally isometric embedding of CAT(0) cubical complexes. We refer the 
reader to the final sentence in the statement of Theorem 1 of \cite{CW}.

\subsection{Some HNN extensions.}
Other hyperbolic groups can be obtained by taking HNN extensions of 
surfaces, and applying the Bestvina-Feighn combination theorem. These are 
two-dimensional but have boundary not homeomorphic to the circle. An explicit 
construction of such an example is given below. These examples necessarily have
local cut points in the boundary. It is known from work of M. Kapovich and B. Kleiner
\cite{KK} that the boundary of a one-ended Gromov hyperbolic group which is 1-dimensional
and has no local cut points is homeomorphic to either the Sierpinski carpet or the 
Menger curve. (see \cite{KapBen}, section 8, for a further discussion). This 
raises the following question which we are so far unable to resolve:

\begin{question}
Do the groups $\DiffD$, $PB_m$, or $\Mod(S)$, for $S$ a closed orientable 
surface, admit quasi-isometric embeddings of 2-dimensional hyperbolic 
groups with boundary homeomorphic to either a Sierpinski carpet or a 
Menger curve? 
\end{question}

\begin{figure}[htb] 
\begin{center}
\input{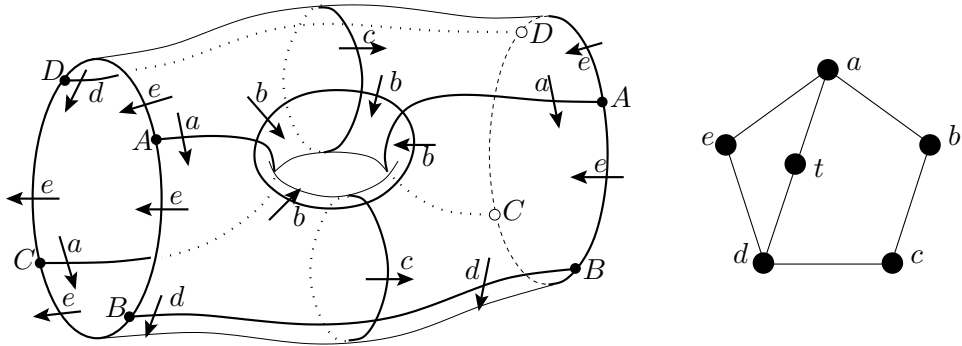}
\end{center}
\caption{We obtain the surface $F$ by identifying the two boundary
circles of the surface on the left. The fundamental group of the
complex $X$ embeds in $G(\Delta)$, where $\Delta$ is the graph on
the right.}\label{F:extension}
\end{figure}

\paragraph{\bf  Construction of HNN extension examples.}
Consider the genus 2 orientable surface $F$ with 
dissection obtained as shown in Figure \ref{F:extension}. This defines a 
quasi-isometric embedding of the group $\pi_1(F)$ in the right-angled 
Artin group $G(C_5)$ where $C_5$ denotes the five cycle graph. Now form 
a double cover $F_2$ of $F$ by cutting along the $e$-curve and gluing 
two copies of the subsurface shown in Figure \ref{F:extension}. Further 
modify the dissection by doubling each of the two $e$-curves (replacing
it with two parallel curves labelled in the same way). 

Now between each pair of parallel $e$-curves we can define a simple closed 
path, labelled $\ga$ and $\ga'$, respectively. We construct a complex $X$
by attaching an annulus $A$ to our surface, joining the two curves $\ga$ 
and $\ga'$. We extend the dissecting curves over the new annulus $A$, by
drawing four segments which connect the two boundary components of $A$,
and which are transversely labelled $a$ and $d$ in the obvious way.
We also introduce a new dissecting curve, namely the core curve of
the annulus, which shall be transverely labelled $t$.

We also modify the right-angled Artin group by adding
a generator $t$ to $G$ which commutes with $a$ and $d$, but not 
with $e$, $b$, or $c$. Thus $G=G(\Delta)$ for a new graph $\Delta$ 
containing $C_5$. It is easy to see that $G(\Delta)$ is of planar type 
(in fact $\Delta^{\rm op}$ is a planar graph).    

Now, the dissection of the surface, extended to the complex $X$, defines a 
homomorphism $\pi_1(X)\to G(\Delta)$. Moreover, by the technique of 
\cite{CW} this can easily be seen to be a quasi-isometric monomorphism: 
the complex $X$ admits a CAT(0) squaring $X_Q$ dual to the dissection, and 
the obvious labelling on the edges of this squaring determines a locally 
isometric embedding of $X_Q$ into the standard cubical complex associated 
to $G(\Delta)$ (see \cite{CW} for more details).

Finally we observe that $\pi_1(X)$ is an HNN-extension. Choosing a basepoint 
in $F_Q^2$ and paths in $F_Q^2$ out to the endpoints of a $t$-edge $E$ of 
the annulus $A$ we define elements $g,g'$ and $t$ in $\pi_1(X)$ 
corresponding to the loops $\ga,\ga'$ and the path $E$, respectively. 
We then have 
\[
\pi_1(X)=\pi_1(F_Q^2)\star_{t:\,\<g\>\to\<g'\>,\,t(g)=g'}=
\pi_1(F_Q^2)\star\<t\>/(tgt^{-1}=g')\,.
\]
It follows from the Bestvina-Feighn Combination Theorem \cite{BeFe}
that $\pi_1(X)$ is a word hyperbolic group. (Note that since the curves 
$\ga$ and $\ga'$ are non-parallel geodesics in $F_Q^2$ the hypotheses of 
\cite{BeFe} Corollary 2.3 are readily satisfied -- namely, no powers of 
$g$ and $g'$ are conjugate in $\pi_1(F_Q^2)$, and $g,g'$ are not proper 
powers of other elements).

\subsection{The commutator subgroup of a right-angled Coxeter group.}

In this subsection we present a natural quasi-isometric embedding of the
commutator subgroup of a right-angled Coxeter group into the corresponding
Artin group. We note that many Coxeter groups are Gromov hyperbolic groups.
In particular, Januszkiewicz and \'Swi\c{a}tkowski \cite{JS} show that 
there exist Gromov hyperbolic right-angled Coxeter groups of virtual 
cohomological dimension $n$, for all $n\geqslant 1$. 

Fix an arbitrary simplicial graph $\Delta$ and let $G(\Delta)$ be the 
associated right-angled Artin group, with standard generators 
$a_1,\ldots,a_n$. One may define the corresponding right-angled Coxeter 
group by adding the further relations that each generator has order 2:
\[
W(\Delta) = G(\Delta)/\<a_i^2:i=1,\ldots,n\>\,.
\]
In the following we shall simply write $W=W(\Delta)$ and $G=G(\Delta)$. 
Observe that $W$ is a group with $2$-torsion. 
However, its commutator subgroup $[W,W]$ has 
index $2^n$ and is torsion free (as a consequence, for instance, of the 
following Lemma). An element $w\in W$ lies in $[W,W]$ precisely when each 
letter $a_i$ appears an even number of times in any word representing $w$. 

\begin{lemma}\label{Coxlemma}
The group $[W,W]$  embeds quasi-isometrically in $G$. 
\end{lemma}

\begin{proof}
We define a map $\phi\co [W,W]\to G$ as follows. Let $w=b_1b_2\ldots b_r$ 
be any word in the usual generators for an element $w\in [W,W]$.
Thus $b_i\in\{ a_1,\ldots,a_n\}$. For each $i\in\{1,\ldots,r\}$
we define $\epsilon_i\in\{\pm 1\}$ by $\epsilon_i=(-1)^{d+1}$ if $b_i$ is 
the $d$th occurrence of that particular letter in the word $w$. Then we set
$\phi(w)=b_1^{\ep_1}b_2^{\ep_2}\ldots b_r^{\ep_r}\in G$. This gives a 
well-defined function since the element $\phi(w)\in G$ is invariant under 
modification of $w$ by trivial insertion or deletion of a square $a_i^2$ 
and substitutions $a_ia_j \leftrightarrow a_ja_i$ for $i,j$ adjacent in 
$\Delta$. Moreover, since each letter $a_j$ appears an even number of 
times in any word for  an element $w\in[W,W]$, the above map defines a 
\emph{homomorphism} $\phi\co[W,W]\to G$. It is clear, since minimal length 
(or reduced) words for elements in $[W,W]$ are mapped to reduced words for 
elements in $G$, that $\phi$ is an injective homomorphism and a 
quasi-isometric embedding.
\end{proof}

\subsection{A closed hyperbolic 3-manifold group.}\label{subsect:3mfld}

Let $\Ups$ denote a regular dodecahedron in $\H^3$ with
dihedral angles $\pi/2$, and define $W_\Ups$ to be the right-angled Coxeter 
group generated by reflections in the faces of $\Ups$. The defining graph 
for $W_\Ups$ as a right-angled Coxeter group shall be denoted $\Delta_\Ups$, 
and is isomorphic to the 1-skeleton of the icosahedron. We observe that 
$\Ga_\Ups:=[W_\Ups,W_\Ups]$ is the fundamental group of a closed oriented
hyperbolic 3-manifold.

\begin{cor}
The closed hyperbolic 3-manifold group $\Ga_\Ups =[W_\Ups,W_\Ups]$ 
may be embedded quasi-isometrically in a right-angled Artin group, 
namely in $G(\Delta_\Ups)$.

\end{cor}
\begin{figure}[htb] 
\begin{center}
\input{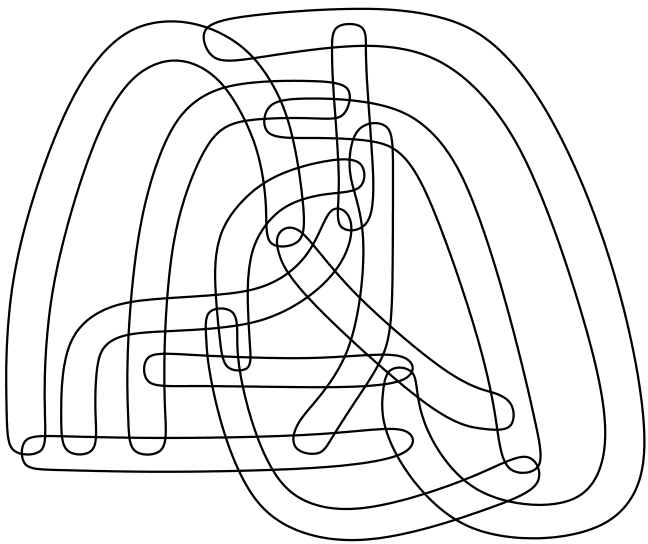}
\end{center}
\caption{A planar circle diagram for the icosahedral graph $\Delta_\Ups$.}
\label{crazycurves}
\end{figure}

Moreover, Figure \ref{crazycurves} shows that the defining graph 
$\Delta_\Ups$ is in fact of planar type. (This would be far from obvious 
at first sight. In fact, while the icosahedral graph $\Delta_\Ups$ is 
planar, its complementary graph is not -- it is obtained from 
$\Delta_\Ups$ by adding an extra edge joining each antipodal pair of 
vertices, and this graph contains a $K_{3,3}$ graph embedded as a subgraph). 
Observe that in Figure \ref{crazycurves} there are 12 circles, 
and each is disjoint from exactly 5 others whose incidence graph is in each 
case a 5-cycle. Thus we have here a circle 
diagram for a regular graph on 12 vertices with vertex valence 5, in which
the link of every vertex spans a 5-cycle (the 5-cycle being self-dual).
In other words, we have a circle diagram for the icosahedral graph $\Delta_\Ups$.

Thus the closed hyperbolic 3-manifold group $\Ga_\Ups$ 
embeds quasi-isometrically in a pure braid group and in the group $\DiffD$.

By Theorem 2 of \cite{JS}, a Gromov hyperbolic right-angled Coxeter group 
can be a virtual $n$-manifold group only if $n\leqslant 4$. 
An example in dimension 4, similar to the $3$-manifold group discussed 
above, would be provided by considering the group of reflections in the 
faces of a hyperbolic hyperdodecahedron (120-cell) with dihedral angles 
$\frac{\pi}{2}$. The commutator subgroup of this reflection group is a 
torsion free subgroup of index $2^{120}$, and the fundamental group of a 
closed hyperbolic $4$-manifold. While this group embeds in a closed surface 
mapping class group (by Corollary \ref{MainCor}(i) and 
Lemma \ref{Coxlemma}), we do not know whether it can be embedded in a braid
group. To obtain such an embedding, it would be sufficient to find a planar 
circle diagram for the corresponding right-angled Artin group on 120 generators.

\begin{ackn}
We would like to thank Daryl Cooper for suggesting the example 
of Subsection \ref{subsect:3mfld}.
We also thank Benson Farb for pointing out to us the work of Januszkiewicz 
and \'Swi\c{a}tkowski \cite{JS} and its relevance to the present work, and in 
particular the existence of hyperbolic 4-manifold subgroups of mapping class groups. 
J.C.~would also like to thank Benson Farb for a number of discussions 
on surface subgroups in mapping class groups which motivated, in particular,
the statement of Corollary \ref{MainCor} (ii).
\end{ackn}

\Addresses

\begin{thebibliography}{999}

\bibitem{Baud}
{\bf A Baudisch,}
{\it Kommutationsgleichungen in semifreien Gruppen}, 
Acta Math. Acad. Sci. Hungar. 29 (1977) 235--249

\bibitem{BG}
{\bf M Benaim, J-M Gambaudo,}
{\it Metric properties of the group of area preserving diffeomorphisms}, 
Trans. Amer. Math. Soc. 353 (2001), no. 11, 4661--4672

\bibitem{BeFe}
{\bf M. Bestvina, M. Feighn,}
Addendum and correction to: ``A combination theorem for negatively 
curved groups'' [\emph{J. Differential Geom.}  35 (1992), no. 1, 85--101]. 
\emph{J. Differential Geom.} 43 (1996), no. 4, 783--788

\bibitem{CW}
{\bf J Crisp, B Wiest,} 
{\it Embeddings of graph braid and surface
groups in right-angled Artin groups and braid groups},
Algebr.~Geom.~Toplogy 4 (2004) 439--472

\bibitem{Dehn}
{\bf M Dehn}, {\it Die Gruppe der Abbildungsklassen},
Acta Math.~69 (1938), 135--206

\bibitem{DynnWie}
{\bf I Dynnikov, B Wiest}, {\it On the complexity of braids}, preprint
{\tt arXiv:math.GT/0403177}.

\bibitem{JS}
{\bf T Januszkiewicz, J {\'S}wi{\c{a}}tkowski},
{\it Hyperbolic {C}oxeter groups of large dimension},
Comment.\ Math.\ Helv. 78 (2003), 555--583

\bibitem{KK}
{\bf I Kapovich, N Benakli,} {\it Hyperbolic groups with low-dimensional boundary},
Annales Sci. Ecole Normale Sup. 33 (2000) No.5, 647--669

\bibitem{KapBen}
{\bf M Kapovich, B Kleiner,} {\it Boundaries of hyperbolic groups},
Combinatorial and geometric group theory 
(New York, 2000/Hoboken, NJ, 2001),  39--93, \emph{Contemp. Math.} 296

\bibitem{Morita}
{\bf S Morita,}
{\it Characteristic classes of surface bundles},
Bull.\ A.\ M.\ S.\ 11, no.\ 2 (1984), 386--388

\end{thebibliography}
\end{document}